**Solving Resource Recharging Station Location-routing Problem through a Resource-space-time Network Representation**

Gongyuan Lu[1,2], Xuesong Zhou[3]*, Qiyuan Peng[1,2], Bisheng He[1], Monirehalsadat Mahmoudi[3], Jun Zhao[1]

*Abstract*: The resource recharging station location routing problem is a generalization of the location routing problem with sophisticated and critical resource consumption and recharging constraints. Based on a representation of discretized acyclic resource-space-time networks, we propose a generic formulation to optimize dynamic infrastructure location and routes decisions. The proposed integer linear programming formulation could greatly simplify the modeling representation of time window, resource change, and sub-tour constraints through a well-structured multi-dimensional network. An approximation solution framework based on the Lagrangian relaxation is developed to decompose the problem to a knapsack sub-problem for selecting recharging stations and a vehicle routing sub-problem in a space-time network. Both sub-problems can be solved through dynamic programming algorithms to obtain optimal solution. A number of experiments are used to demonstrate the Lagrangian multiplier adjustment-based location routing decision making, as well as the effectiveness of the developed algorithm in large-scale networks.

---

[1] School of Transportation and Logistics, Southwest Jiaotong University, Chengdu, Sichuan, China, 610031.
[2] National United Engineering Laboratory of Integrated and Intelligent Transportation, Chengdu, Sichuan, China, 610031.
[3] School of Sustainable Engineering and the Built Environment, Arizona State University, Tempe, AZ, 85281,
* Corresponding author, School of Sustainable Engineering and the Built Environment Ira A. Fulton Schools of Engineering, Arizona State University, Tempe, AZ, 85281, xzhou74@asu.edu, Tel: +1-480-965-5827

## 1. Introduction

### 1.1 Overview

To provide high quality transportation services for customers considering a reasonable cost and budget, transportation planning agencies and logistics companies need to design, operate, and maintain an effectively integrated transportation service and infrastructure network. By optimizing vehicle routes, service frequency, and timetables, the task of transportation service network design aims to satisfy time-dependent customer origin-destination demand subject to various forms of complex resource constraints (in terms of battery capacity, fuel, and working hour durations for drivers and crew). As a result, different types of vehicle service infrastructure locations such as Electric Vehicle (EV) battery recharging stations, locomotive refueling terminals, and railcar inspection depots should be strategically located in a transportation network to meet the required vehicle service requirements.

To optimize the infrastructure location and vehicle routes simultaneously, the Location Routing Problem (LRP) has been extensively studied. Interested readers are referred to the survey papers by Min et al. (1998) and Nagy and Salhi (2007) for the related taxonomy and classification. In this research, we are particularly interested in the joint resource-recharging station location and routing problem (RRS-LRP) with sophisticated constraints on resource consumption and recharging. This problem covers the classical LRP as a special case, where a depot can be viewed as a special type of resource recharging stations (RRS) which provides the needed platform for vehicles to start and end their tours. From a much broader modeling perspective, not only we should consider different vehicle types, namely passenger cars, buses, trucks, aircrafts, and locomotives, but also various forms of resources consumed by vehicles continuously while travelling; to name a few, energy resources including petroleum, diesel, and electricity; driver time resources expressed as working hour duration; as well as driving distance resources associated with regular technical inspection requirements for rail cars or commercial fleet. Accordingly, the generic formulation to be developed hopes to cover various forms of RSS, e.g., battery swapping station, electricity recharging station, gas refilling station, technical maintenance depots, and even hotels for crew rest.

**1.2 Emerging applications of joint resource recharging location-routing problem**

Since a wide range of multi-modal transportation applications can be casted as a vehicle routing problem with resource constraints, our research on the joint RRS-LRP optimization is motivated by two representative and emerging applications. First, we consider EV routing and recharging infrastructure planning in an urban traffic network setting, followed by locomotive routing and refueling station optimization problem for regional rail infrastructure network design.

The continuous advance and development of electric vehicles hold the promise of meeting the daily urban travel requirements while offering a promising way to reduce local traffic emissions and petroleum dependence. As recognized by many planners and researchers (to name a few, Mak et al. 2013 and He et al. 2013), the wide use of plug-in hybrid electric vehicles and all-electric vehicles requires systematic infrastructure network planning and sufficient deployment of EV charging stations. Two types of charging strategies are now commonly available including plug-in charging and battery swapping. Specifically, battery swapping stations (BSS) enable a large number of EVs to be recharged very fast during peak hours, while plug-in charging stations can provide a cost-effective solution if much longer time duration (e.g. 15 minutes to 10 hours) is allowed. Given different cost effectiveness trade-offs, the types of charging equipment and the location of charging stations have to be carefully selected and optimized to meet the critical needs of urban drivers and fleets in a cost-effective fashion.

How to select the type of fuel supply and the locations of refueling stations is one of strategical decisions faced by railroad companies around the world. For instance, North American railroads spent approximately $11.6 billion on 4.1 billion gallons of diesel in 2013, and significant attention is dedicated to how to improve efficiency and deploy various fuel cost reduction strategies (Gladstein, Neandross & Associates, 2014). Railroad companies to decide where and how many times to refuel the locomotives to minimize the overall costs due to location-dependent fuel price, examined by Nourbakhsh and Ouyang (2010). Similarly, in a high-speed rail network, after traveling a certain distance (let say 4,000 km), each train unit must travel back to a depot for performing mandatory maintenance and inspection activities. In this case, the maximum travel distance becomes a hard resource constraint which is "a train has only 4,000 km distance-resource to consume, and it must be recharged at one of the inspection depots when the remaining distance-equivalent resource is insufficient". Faced with tight investment budget constraints, high-level decisions have to be made to carefully locate train unit maintenance depots, while providing sufficient service coverages for time-sensitive and spatially distributed passenger transporting demand.

**1.3  Review of related decomposed sub-problems and solution algorithms: math details**

The RRS-LRP integrates a number of classical optimization problems, such as the location problem and the vehicle routing problem with an embedded Resource Constrained Shortest Path (RCSP) sub-problem as its key building block. The generic RCSP has been well studied by Handler and Zang (1980), and common solution approaches include Lagrangian relaxation, path ranking-based approach, dynamic programming strategies, and branch and bound method examined by Feillet et al. (2004); Santos et al. (2007); Carlyle et al. (2008); Pugliese and Guerriero (2013a). Focusing on precisely tracking resource consumption states along the traveling path, Feillet et al. (2004) proposed a label correcting-based exact solution procedure to solve the elementary shortest path problem. Multiple types of resources are involved in the problem studied by Bektas and Laporte (2011) and Pugliese and Guerriero (2013b).

The shortest path problem with resource constraint (SPPRC) was systematically examined in the modeling framework proposed by Desrochers et al. (1986), where time, car load, and break duration can be considered as resources which varies along a path according to a resource extension function. They also provided a systematic taxonomy for SPPRC, classified by resource accumulating process, path-structural constraints, objective, and underlying network. Typically, cycles need to be eliminated in a physical network while solving the SPPRC by formulating the problem over acyclic 2-dimensional time-space networks can eliminate cycles naturally.

Compared to RCSP and SPPRC, the vehicle routing problem with recharging station (VRP-RS) is more complex in its own right due to additional dimensions of multiple vehicles and time-sensitive demand satisfaction requirements. To capture the resource consumption and recharging dynamics, a set of linear constraints are typically needed to track resource usages on each transportation link along the path (Erdogan and Miller-Hooks, 2012; Hiermann et al. 2014; He et al. 2014; Wang et al., 2014, Schneider et al., 2014. The recent study by Hiermann et al. (2014) and Schneider et al. (2014) takes demand time windows into account, and Schneider et al. (2014) further considered the different charging time durations at different battery levels.

As mentioned before, there are a wide range of resource types tightly connected with the vehicle routing problem; for example, locomotive fuel (Nourbakhsh and Ouyang, 2010), normal vehicle fuel (Conrad and Figliozzi, 2011), alternative vehicle fuel (Erdogan and Miller-Hooks, 2012), distance (Berger et al. 2007), crew duty (Steinzen et al. 2010) and electricity power (He et al. 2014; Worley et al. 2012; Schneider et al. 2014; Yang and Sun, 2015). The objective functions in the context of RRS-LRP usually covers the total transportation cost and fixed RSS construction cost (Berger et al. 2007; Worley et al. 2012; Yang and Sun. 2015), as well as crew operating cost (Steinzen et al. 2010) and number of vehicles (Conrad and Figliozzi, 2011).

To the best of our knowledge, very limited existing literature clearly defines and solves RRS-LRP, but there are a variety of related studies involving resource recharging location and routes decisions. To provide a systematic comparison of key modeling components, we bring Table 1 to examine various resource and vehicle definitions, in conjunction with model formulations and solution algorithms. Specifically, Berger et al. (2007) proposed a model of LRP with distance constraints and given feasible route set connected to a depot. In a systematic study by Nourbakhsh and Ouyang (2010) based on real world railroad applications, a mixed integer programming model is developed to determine optimal locations of contracted fuel stations based on pre-determined locomotive routing profiles. A time-space network based vehicle and crew scheduling model is proposed by Steinzen et al. (2010), where the duties generation is modeled as a resource constrained shortest path problem, and crew must break for enough time at relief points. The model is solved by a combined column generation and Lagrangian relaxation algorithm. Worley et al. (2012) studied a simultaneous vehicle routing and charging station siting problem, in which routes between two charging stations are defined as "parts" connecting the customers with traveling length limitations. Mak et al. (2013) considered an infrastructure planning problem for EV battery swapping stations with both RSS construction budget constraint and recharging requirements along EV routes. In the battery swapping station location-routing problem examined by Yang and Sun (2015), the resource tracking variables are adopted to represent battery states for a vehicle entering / leaving a node.

This paper aims to formulate a general class of the resource constrained location routing problem for minimizing the transportation cost with given time-dependent demand constraints and recharging station capacity constraints. To capture energy consumptions as a function of vehicle driving speed, a number of recent studies (Yang and Zhou,2014) explicitly consider a time-expanded network where the transportation network is replicated in discrete time intervals. Mahmoudi and Zhou (2015) develop a new state-space-time hyper-network representation to solve the vehicle routing problem with pickup and delivery services and time windows. In addition to the topological network at each time interval, they introduced a vehicle-specific state for indicating the individual passengers that are on board, and such a construct enables them to use a forward-pass dynamic program for solving the decomposed vehicle routing sub-problems.

In this research, through an appropriate use of additional dimension of resource, , we construct a directed acyclic resource-space-time (RST) network to reformulate RRS-LRP as a multi-commodity flow problem with linking constraints to ensure that all transportation demands and recharging requirements are satisfied. This new RST representation offers a number of contributions to address a number of critical issues in vehicle routing with hard resource recharging constraints. First, it can explicitly embed or pre-build resource consumption and recharging constraints into a well-structured formulation, so that the vehicle trajectories in the hyper network automatically reflect time-sensitive resource changes in both transporting and recharging processes. Second, by using the recharging activity and demand on arcs to consider the location-related activities, we are able to cast the location-routing problem as seamlessly integrated model with a time-dependent network flow problem for vehicle routing and a conventional knapsack problem for selecting recharging station locations with limited budget. Lastly, these two sub-problems are solved within a computationally efficient Lagrangian decomposition framework through iteratively adjusting two sets of Lagrangian multipliers for time-dependent demand links and recharging station links.

The remainder of this paper is organized as fo1llow. In section 2, we provide a problem description within a resource-space-time network construct. The model of RRS-LRP is formulated in section 3 to capture all essential constraints for vehicle routing and recharging station selection. Section 4 presents a Lagrangian decomposition algorithmic framework to solve two sub-problems by using dynamic programming. In section 5, we systematically discuss the differences between our proposed model and other alternative formulations, including the rudimentary model and pre-determined route model for RRS-LRP. Section 6 conducts numerical experiments to examine the effectiveness of the proposed formulation and algorithm.

Table 1. Summary of related modeling and problem solving methodology

| Publication | Type of Problem | Type of Vehicle | Resource recharging station / depot | Routing resources, and business rules | Demand representation | Objective function | Solution algorithm |
|---|---|---|---|---|---|---|---|
| Berger, Coullar, and Daskin, 2007 | LRP with distance constraint | Generic vehicle | N/A | Distance, implicit recharging at same facility | On node | Minimize total facility cost and transportation cost | BP |
| Nourbakhsh and Ouyang, 2010 | Locomotive refueling station location and routing | Locomotive | Fuel station with limited capacity | Diesel refueling on node | Pre-determined locomotive trip | Minimize total fuel cost, delay cost and fixed contract cost of fuel station | LR |
| Steinzen et al., 2010 | Vehicle- and Crew-scheduling problem | Vehicle | Relief point of crew | Cost of crew duty | Demand on link | Minimize summation of vehicle cost and crew cost | LR&CG |
| Conrad and Figliozzi. 2011. | Recharging VRP | Vehicle | Customer location | Fuel | Demand on node | Primary objective: number of vehicles Secondary objective: total cost | H |
| Worley et al. 2012 | EVs Recharging station location routing problem | Electric vehicle | Recharging station | Electricity, recharge on node | Demand on node | Minimize total transportation, recharging, and charging station placement costs | N/A |
| Erdogan and Miller-Hooks. 2012 | Green vehicle routing problem | Alternative fuel vehicle | Alternative fuel station | Alternative fuel, refuel on node | Demand on node | Minimize total travel distance | H |
| Schneider et al. 2014 | Electric VRP with time window and recharging station | EV | Recharging station | Battery, refuel at node | Demand on node | Minimize total travel distance | H |
| Yang and Sun, 2015 | EVs Battery swap station LRP | Electric vehicle | Battery swapping station | Battery, swapping on node | Demand on node | Minimize construction cost and EV shipping cost | H |

| Our paper | Resource recharging station LRP | Generic vehicle | Resource recharging station | Recharging on link | Demand on link | Minimize total transportation cost | LR&DP |

Solution algorithm: LR – Lagrangian Relaxation, H – Heuristic method, CG – Column Generation, BP – Branch and Price, DP-Dynamic Programming

## 2. Problem statement with a resource-space-time network representation

We now start formally defining the RRS-LRP as follow. Given candidate resource recharging station locations and transportation demand in terms of origins, destinations, departure time and expecting arrival time, the problem studied in this paper aims to determine the joint decision of RRS locations and vehicle routes with a goal of minimizing the total transportation cost, subject to the resource recharging requirements along vehicle routes, customer demand satisfaction constraints, and RRS recharging volume capacity and construction budgetary constraints. The notations used in the RRS-LRP model are first listed in Tables 2-4.

Table 2. Sets and indexes for RRS-LRP model

| Symbol | Definition |
|---|---|
| $N$ | Set of resource-space-time nodes |
| $E$ | Set of resource-space-time links |
| $N_s$ | Set of nodes of resource recharging stations in RST network |
| $E_s$ | Set of links that connect the resource recharging station nodes in RST network |
| $N^p$ | Set of physical nodes |
| $E^p$ | Set of physical links |
| $N_r^p$ | Set of physical recharging station nodes |
| $E_r^p$ | Set of physical links connect the recharging station |
| $R$ | Set of resource indexes in resource dimension |
| $S$ | Set of physical node indexes in space dimension |
| $T$ | Set of time indexes in time dimension |
| $V$ | Set of vehicles |
| $\Psi$ | Set of demands in space-time network |
| $i,j$ | Index of physical nodes |
| $k$ | Index of resource recharging station nodes |
| $t,t'$ | Index of time intervals |

| Symbol | Definition |
| --- | --- |
| $r, r'$ | Index of resource intervals |
| $v$ | Index of vehicles |

Table 3. Given parameters in the optimization model

| Symbol | Definition |
| --- | --- |
| $\varphi(i,j,t,t')$ | Transportation service demand flow from node *i* at time *t* to node *j* at time *t'*, =1 when the demand exists on space-time link $(i,j,t,t')$; =0 otherwise. |
| $c_{i,j}$ | Travel cost of link $(i,j)$ |
| $r_{i,j}$ | Cost of resource if travel from node i to j |
| $c_{i,j,t,t'}$ | Travel cost of space-time link $(i,j,t,t')$ |
| $r_{i,j,t,t'}$ | Consumed resource if travel through space-time link $(i,j,t,t')$ |
| $c_{i,j,t,t',r,r'}$ | Cost on resource-space-time link $(i,j,t,t',r,r')$ |
| $p_k$ | Maximum capacity of recharging station *k* |
| $q_v$ | Maximum amount of resource carried by vehicle *v* |
| $TT_{i,j}, TT_{i,j}^{min}, TT_{i,j}^{max}$ | Travel time, minimum travel time and maximum travel time on link $(i,j)$ |
| $TT_{i,j,r,r'}$ | travel time from node $i$ at resource status $r$ to node $j$ at resource status $r'$ |
| $o(v), d(v)$ | Origin node and destination node of vehicle *v* |
| $r_v$ | Origin resource status of vehicle *v* |
| $r_0$ | Minimum resource status |
| $T_v, T'_v$ | Designated departing time and arrival time of vehicle *v* |
| $b_k$ | Construction cost of depot *k* |
| $B$ | Total depot construction budget |

Table 4. Decision variables to be optimized

| Symbol | Definition |
| --- | --- |
| $x^v_{i,j,t,t',r,r'}$ | Binary vehicle routing variable; =1, if vehicle *v* travel on RST link $(i,j,t,t',r,r')$; =0, otherwise. |
| $w_k$ | Binary depot location decision variable, =1, if resource recharging station k is selected to be established; =0, otherwise. |

Before elaborating the RRS-LRP formulation, we initially describe the Resource-Space-Time (RST) network, which is a combination of three dimensions of resource states, nodes and links in the space layer, and the time transition layer. Let $RSTG = (N, E)$ be a directed resource-space-time network with *N* as the RST nodes set and *E* as the RST link set. To construct this hyper network, the node *i* in the physical network is extended to a RST node $(i, t, r)$, indicating that the vehicle maintains a resource level of *r* at physical node *i* at time *t*. Accordingly, a RST link with a form of $(i, j, t, t', r, r')$ corresponds to a vehicle traveling or resource recharging activity through physical link $(i,j)$ by taking a time duration $(t'-t)$ and a resource change $(r'-r)$. The constructed RST network assumes all activities to be performed on the links, this activity-on-link modeling framework requires us to transform the physical node, which represents the recharging station to a link $(i,j)$, as the network transformation illustrated in Fig. 1. Specifically, the activity performing time $(t'-t)$ in the network is always positive for either traveling or recharging activities, while the change of resource level $(r'-r)$ is negative on transporting links and positive on recharging links. The features of such an activity-on-arc representation is examined in details by Assad and Golden (1995). In this paper, we also define the demand of customers on links.

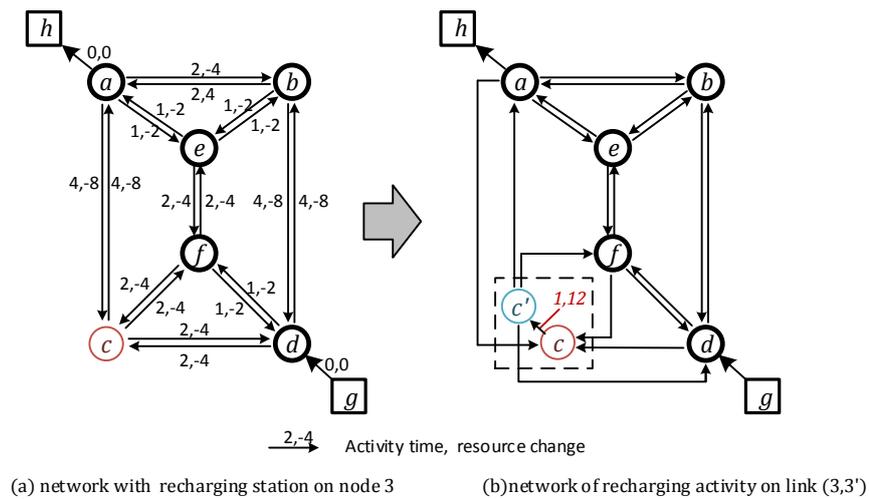

(a) network with recharging station on node 3   (b) network of recharging activity on link (3,3')

Figure 1. Transforming the recharging node to recharging link

Fig 2. illustrates a RST network on the physical network with 3 nodes and 2 links, as well as a maximum level of the resource storage in a vehicle as 4 units. Figs 2(b) and 2(c) show the network projections to the 2-dimensional space-time (S-T) and resource-time (R-T) places, respectively. The RST route, shown in Fig. 2(d), starts at time 1 and ends at time 4 in the 3-D RST network, corresponding to a sequence of RST nodes $(a,1,3) \to (b,3,1) \to (c,4,4)$. The multi-

dimensional RST link structure $(i,j,t,t',r,r')$ explicitly codes the traveling / recharging time duration and resource changes, while the node label $(i,t,r)$ can be directly used to trace the cumulative time spent and resource storage level for each vehicle.

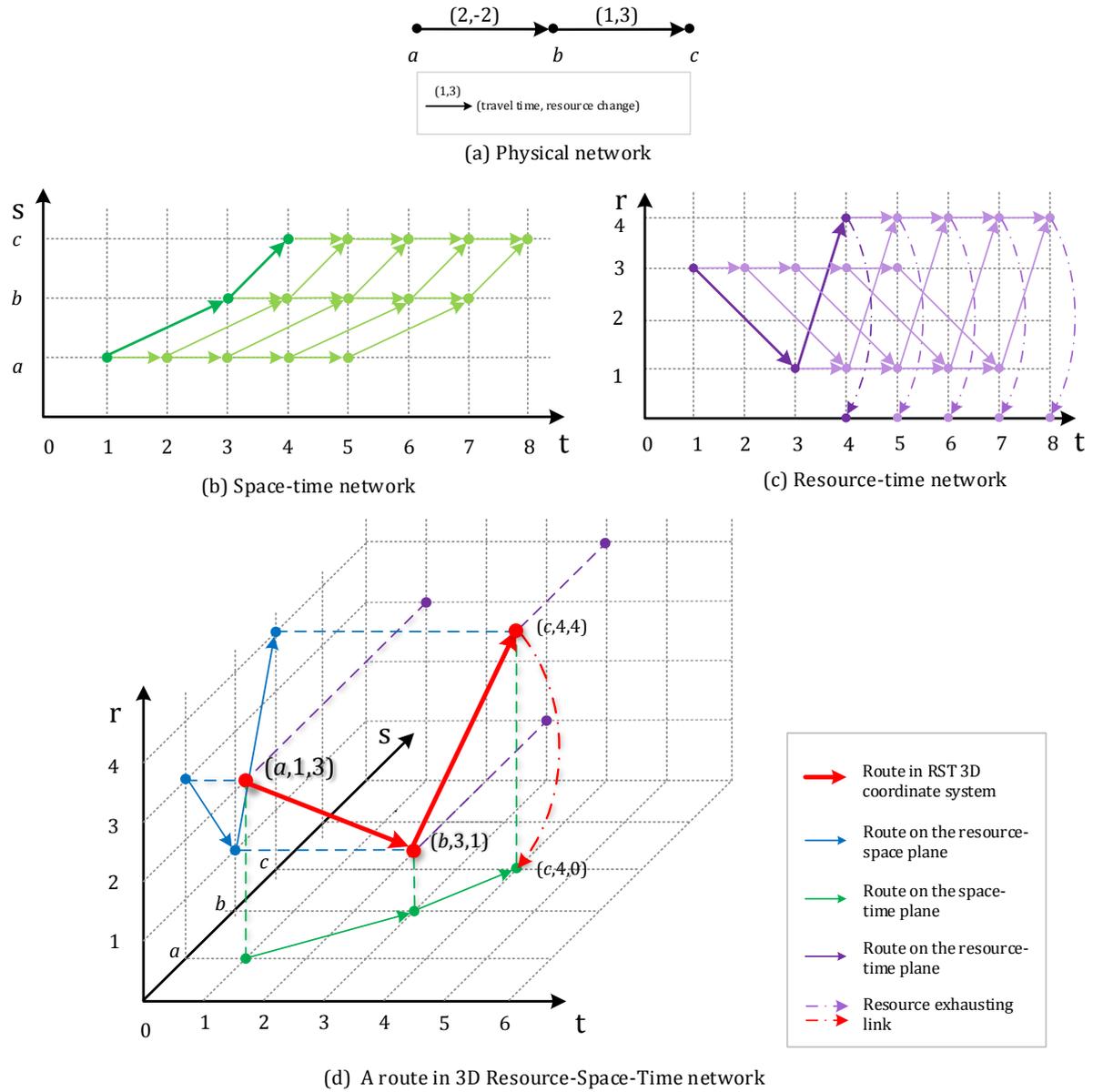

Figure 2. A simple example of resource-space-time network

When recharging EVs, the increase of the battery volume is dependent on the total recharging time and starting battery level. Fig. 3 shows a discretized version of the battery volume change as a function of recharging time duration, for different types of recharging equipment. Obviously, a vehicle does not need to enter a recharging station only if the fuel tank is empty, so the resource recharging link can start with different initial resource states as shown in Fig. 3

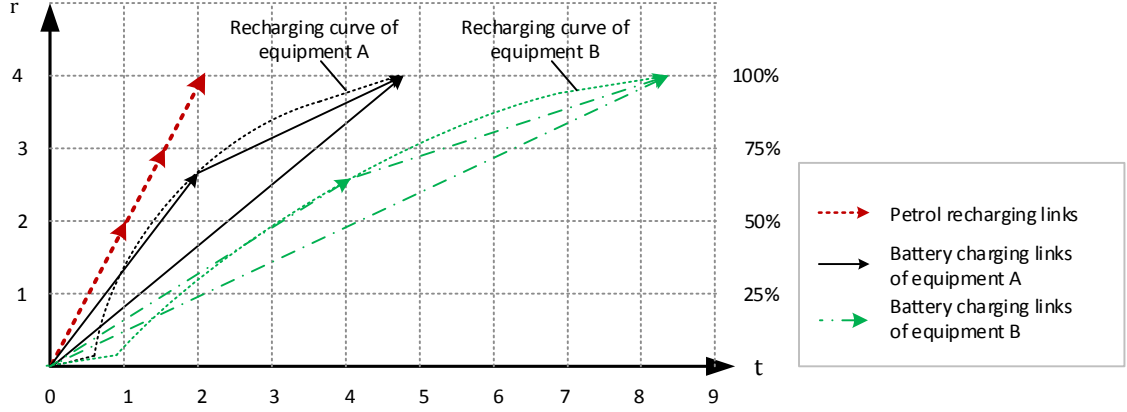

Figure 3. Recharging links with different types of equipment

Similarly, since the energy or resource consumption in general depends on travel time or speed, one can easily represent a discretized version of the energy consumption rate as a function of EV's speed *V* per time interval according to Eq. (1).

$$\Omega = -0.064 + 0.0056V + 0.00026(V-50)^2 \qquad (1)$$

While traveling on the same route, the vehicle can either accelerate to arrive the destination earlier or decelerate to arrive later. The resource consumption rate changes associate with the speed variation, have been shown in Fig. 4, at speed of 40, 60, and 80 km/h. The vehicle finishes the same route with different time and resource consumption. The remaining resource status is labeled next to the space-time node. In the RST network, we can use a RST link to demonstrate the time and resource consumption simultaneously. For instance, in the vehicle route of speed 80 km/h, the space-time link from node (*d*, 5) at resource status 12 to node (*a*, 9) at resource status 3 can be represented by a single RST link (*d*,5,12)→(*a*,9,3).

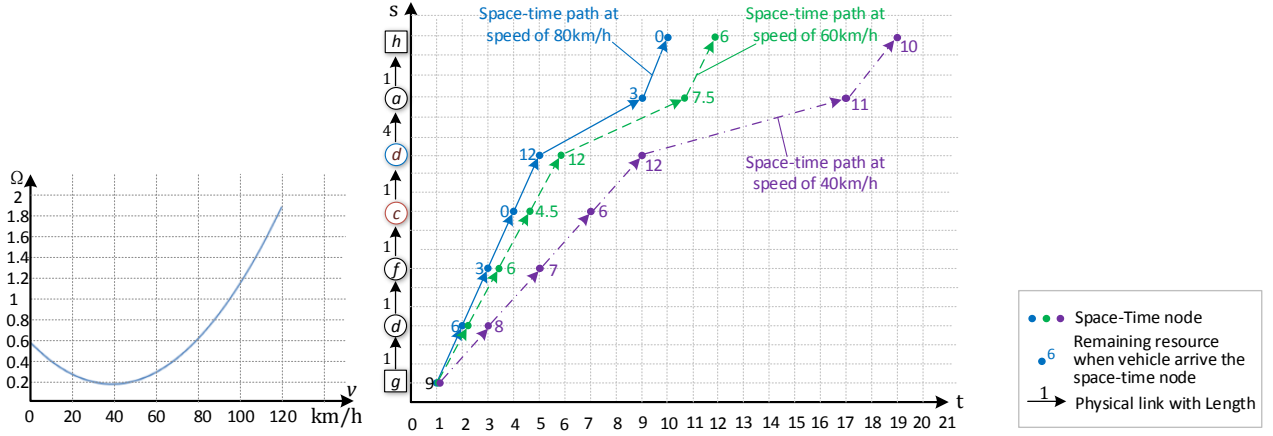

(a) Resource consumption process with respect to the CO emission

(b) Resource change processes at different speed via the same route

Figure 4. Time sensitive resource consumption process in RST network

## 3. Model of Resource Recharging Station Location Routing Problem(RRS-LRP)

By defining binary vehicle routing variable $x^v_{i,j,t,t',r,r'} = 1$, where vehicle $v$ travels on RST link $(i,j,t,t',r,r')$, we now proceed to the mathematical programming formulation. The objective function (2) aims to minimize the total traveling costs, subject to the flow balance constraints, the RRS service capacity constraints, the customer demand constraints, and the construction budget constraints. Without loss of generality, we assume that the vehicles start and end their routes at any nodes during the planning horizon, and can be recharged at any resource recharging station. Using the algorithm in Appendix, a complete RST network can be built based on the physical network.

**Model 1. RRS-LRP model based on RST network**

**Objective Function** $z = \min \sum_{v \in V} \sum_{(i,j,t,t',r,r') \in E} (c_{i,j,t,t'} \times x^v_{i,j,t,t',r,r'})$  (2)

**Subject to.**

**(1) Flow balance constraints**

For vehicle $v$, the flow balance constraint at RST node ($j$, $t'$, $r'$) can be formulated as:

$$\sum_{(j,t',r')\in N} x^v_{i,j,t,t',r,r'} = 1, \forall v \in V, i = o(v), t = T_v, r = r_v$$

(3)

Eq. (3) ensures that only one vehicle leaves the given RST origin node $(i,t,r)$ of vehicle v

$$\sum_{r'\in R}\sum_{(i,t,r)\in N} x^v_{i,j,t,t',r,r'} = 1, \forall v \in V, j = d(v), t' = T'_v, r' = r_0$$

(4)

Eq. (4) ensures that only one vehicle with a specific arrival time and resource state arrives to the physical destination node. Notice that the remaining resource status after a vehicle arriving its destination is always unknown, thus the flow balance at the end of a trip needs to be maintained through a set of virtual exhausting arcs for each space-time destination node to allow the vehicles consuming all the remaining resource and return back to the super sink.

$$\sum_{(j,t',r')\in N} x^v_{i,j,t,t',r,r'} - \sum_{(j,t',r')\in N} x^v_{j,i,t',t,r',r} = 0,$$

$$\forall i,t,r \in N, (i,t,r) \notin \{(o(v),T_v,r_v),(d(v),T'_v,r_0)\} \forall v \in V \quad (5)$$

Eq. (5) guarantees the flow balance on non-origin and non-destination RST nodes.

**(2) Resource recharging station capacity constraints**

This constraints state that every recharging station $k$ has its maximum capacity to serve vehicles during the planning horizon.

$$\sum_{v\in V}\sum_{t\in T, r\in R}\sum_{(j,t',r')\in N} x^v_{k,j,t,t',r,r'} \leq p_k \times w_k, \forall k \in N_s$$

(6)

**(3) Demand satisfaction constraints**

The essential purpose of transportation service is to serve time-dependent travel demands, defined through the space-time demand link coefficient $\varphi(i,j,t,t')$.

$$\sum_{v \in V} \sum_{r,r' \in R} x_{i,j,t,t',r,r'}^{v} \geq \varphi(i,j,t,t'), \forall \varphi(i,j,t,t') \in \Psi$$

(7)

**(4) Station construction constraints**

$$\sum_{k \in N_r} b_k \times w_k \leq B$$

(8)

There are also binary definitional constraints for variables $x$ and $w$.

$$x_{i,j,t,t',r,r'}^{v} \in \{0,1\}, \forall v \in V, \forall i,j,t,t',r,r' \in E$$

$$w_k \in \{0,1\}, \forall k \in N_s$$

With the RST network, the RRS-LRP is formulated as a classic LRP formulation system with the resource-space-time variables. The model is implemented in GAMS 24.2.3 to verify its correctness. The above formulation involves $|V||N|^2 + |N_r^p|$ binary variables and $2 + 2|V| + |V|(|S| - 2)|T||R| + |N| + |\Psi|$ equations. For real world problems, the scale of the formulation is very large, so an efficient solution approach is required. Fortunately, the well-organized network and formulation system gives the possibility of developing a customized solution framework to solve the problem efficiently. Virtual vehicle sets are used to account for possible solution infeasibility, the method of applying the virtual vehicle is referred to Mahmoudi and Zhou (2015).

## 4. Solution Approach

### 4.1 Problem Decomposition and Lagrangian Relaxation

In a Location-routing problem (LRP), the location of depots and routes of vehicles must be determined simultaneously. There are a wide range of sequential solution algorithms that usually divide the problem into 2 or 3 stages (Hansen et al., 1994, Lin et al. 2002, etc.). In our research, we plan to apply a Lagrangian decomposition method to decouple the problem into two sub-problems, the resource recharging location problem and resource constrained vehicle routing problem.

Table 5. Formulation structure of constraints

| Group Index | Constraints | Variable involved | Form |
|---|---|---|---|
| 1 | Flow balance constraints | $x^v_{i,j,t,t',r,r'}$ | $AX = B$ |
| 2 | RRS capacity constraints | $x^v_{i,j,t,t',r,r'}, w_k$ | $AX \leq CW$ |
| 3 | Demand link constraints | $x^v_{i,j,t,t',r,r'}$ | $AX \geq B$ |
| 4 | RRS construction budget constraints | $w_k$ | $AW \leq B$ |

Four groups of linear constraints are listed in Table 5 to better illustrate the model structure. Specifically, groups 1 and 3 are route selecting constraints for VRP, while group 4 is the construction budget constraints related to RRS decision variables. The recharging station capacity constraints (group 2) are coupling constraints linking two different sub-problems. In particular, as shown in Fig. (9), we introduce a set of non-negative Lagrangian multiplier $\theta_k$ to dualize the depot capacity constraints of Eq. (6) and introduce a set of non-negative Lagrangian multiplier $\varepsilon_{i,j,t,t'}$ to dualize the demand link constraints of Eq. (7).

Table 6. Notation of Lagrangian multipliers

| Symbol | Definition |
|---|---|
| $\theta_k$ | Lagrangian multiplier corresponding resource recharging station capacity |
| $\varepsilon_{i,j,t,t'}$ | Lagrangian multiplier corresponding demand links |

**Objective Function**

$$L(x,w,\varepsilon,\theta) = \sum_{v \in V} \sum_{(i,j,t,t',r,r') \in E} (c_{i,j,t,t'} \times x^v_{i,j,t,t',r,r'})$$

$$+ \sum_{(i,j,t,t',r,r') \in E} \varepsilon_{i,j,t,t'} \times \left[ \varphi(i,j,t,t') - \sum_{v \in V} x^v_{i,j,t,t',r,r'} \right]$$

$$+ \sum_{k \in N_r} \theta_k \times \left( \sum_{v \in V} \sum_{r \in R, t \in T} \sum_{(j,t',r') \in N} x^v_{k,j,t,t',r,r'} - w_k \times p_k \right) \quad (9)$$

**Subject to:**

Flow balance constraints (3-5) and RRS construction budget constraint (8)

Binary constraints for **x** and **w**, non-negative constraints for **ε** and **θ**

The multiplier $\varepsilon_{i,j,t,t'}$ can be interpreted as the profit of completing the demand satisfying tasks of $\varphi(i,j,t,t')$. With this group of multipliers, the major goal of the objective function (9) is to obtain maximum profit by picking up as many demand links as possible. The multiplier $\theta_k$ is involved in both knapsack sub-problem and vehicle routing sub-problems, which indicates the utility of each depot in the first problem, and can be interpreted as the system marginal cost for consuming capacity resources of depots in the second problem.

The dual problem **Pxw** and its 2 sub-problems are shown below.

**Dual problem Pxw**

$$L(x,w,\varepsilon,\theta) = \sum_{(i,j,t,t',r,r') \in E} \varepsilon_{i,j,t,t'} \times \varphi(i,j,t,t') - \sum_{k \in K} \theta_k \times p_k \times w_k$$

$$+ \sum_{v \in V} \sum_{(i,j,t,t',r,r') \in E} (c_{i,j,t,t'} - \varepsilon_{i,j,t,t'} + \theta_k) \times x^v_{i,j,t,t',r,r'} \quad (10)$$

s. t.    Non-negative constraints for **θ** and **ε**

**Sub-problem SPw: knapsack problem**

$$\max L_{x-2} = \sum_{k \in K} \theta_k \times w_k \times p_k$$

s. t.    RRS construction budget constraints (8)

Binary constraints for **w**

**Sub-problem SPx: Multi-vehicle routing problem with recharging station (VRP-RS)**

$$\min \sum_{v \in V} f_v = \sum_{v \in V} \sum_{(i,j,t,t',r,r') \in E} (c_{i,j,t,t'} - \varepsilon_{i,j,t,t'} + \theta_i) \times x^v_{i,j,t,t',r,r'}$$

s. t.    Flow balance constraints (3-5)    Binary constraints for **x**

In **SPx**, the generalized cost/profit parameters in the brackets can be expressed in terms of $\text{profit}(i,j,t,t')$:

$$\text{profit}(i,j,t,t') = c_{i,j,t,t'} - \varepsilon_{i,j,t,t'} + \theta_i \tag{11}$$

### 4.2 Solution algorithm

A Lagrangian relaxation algorithm framework is plotted in Fig. 5, with a dynamic programming based algorithm for sub-problem SPw. For SPx operated in the context of resource-space-time network, we develop a dynamic programming-based algorithm based on a time-expanded network as an acyclic direct graph.

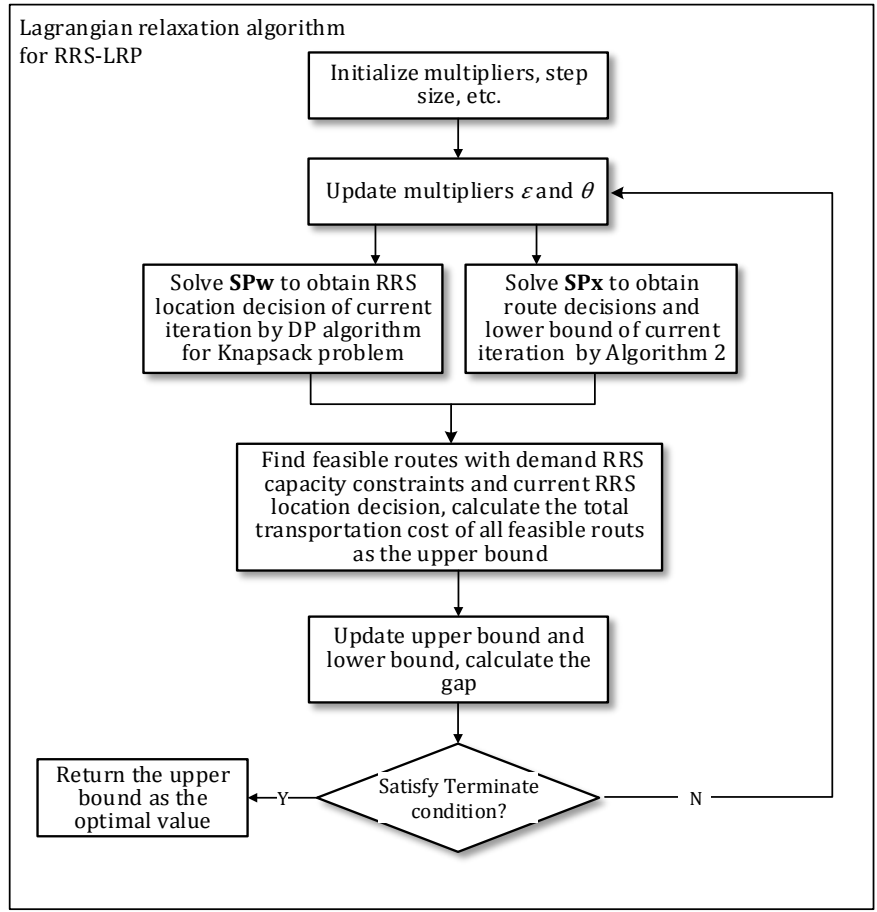

Figure 5. Lagrangian relaxation algorithm framework

**(1)Algorithm 1: Algorithm for single vehicle routing problem with recharging station**

A Dynamic Programming (DP) algorithm is developed based on the directed acyclic RST network. Importantly, the proposed algorithm does not require a network building step, as all the loops automatically scan the outgoing arcs through a forward DP solution approach that allows negative travel cost on links.

To implement the Dynamic Programming algorithm for single VRP-RS, we introduce the following additional notation. Note that, the time-dependent and location-dependent resource consumption/recharging parameter $r_{i,j,t,t+TT_{i,j}}$ is adopted in Algorithm 2, which could be either positive for a resource consumption arc or negative for a resource recharging arc.

Table 7. Additional notation for DP algorithm used to solve VRP-RS

| Symbol | Definition |
|---|---|
| $lc(i,t,r)$ | Label cost at RST node $(i,t,r)$ |
| $lc'(j,t',r')$ | Temporary Label cost at RST node $(j,t',r')$ for recording |
| $(o,t_0,r_0),(d,t_d,r_d)$ | RST node of origin node and destination node in an optimal path |
| $preNode(j,t',r')$ | Predecessors, the index pointed to the previous RST node of $(j,t',r')$ in the optimal path |

**Step 1: Initialization**

For all $lc(i,t,r) = M$

$lc(o,t_0,r_0) = 0$

**Step 2: recursion to find resource constrained optimal path**

For (t =0 to |T|)

    For each exist space-time link $(i,j)$

        For each resource status $r$

            For $(TT_{i,j} = TT_{i,j}^{min}$ to $TT_{i,j}^{max})$

            $r' = r + r_{i,j,t,t+TT_{i,j}}$

            If $(r' > 0)$

                **Label cost updating**

                $lc'(j,t',r') = lc(i,t,r) + profit(i,j,t,t')$ in which $t' = t + TT_{i,j}$

                If $(lc'(j,t',r') \leq lc(j,t',r'))$

                      $lc(j,t',r') = lc'(j,t',r')$

                      $preNode(j,t',r') = (i,t,r)$

                End if

            End if

            End For of $TT_{i,j}$

        End For of $r$

End for of link $(i,j)$

End for of $t$

**Step 3: Track back for the optimal path**

Track back from super sink node $(d, t_d, r_d)$ back to super source node $(o, t_0, r_0)$ using the predecessor sequence $preNode(j, t', r')$.

(2) **Algorithm 2: Lagrangian relaxation algorithm (LRA) for Px and Py**

**Step 1: Initialization**

Start algorithm at iteration $m=1$, initialize demand profit multipliers $\varepsilon_{i,j,t,t'}^m = 0$ and capacity price multiplier $\theta_k^m = 0$. Set the step size $\alpha^m = 1$.

**Step 2: LR Multiplier updating**

2.1: update step size $\alpha^m = 1/(m + 1)$;

2.2: update $\varepsilon_{i,j,t,t'}^m$ with current iteration index and current max profit path solution by sub-gradient method:

For each vehicle $v$

    For each link $(i, j, t, t')$

$$\varepsilon_{i,j,t,t'}^m = \max\left\{0, \varepsilon_{i,j,t,t'}^{m-1} + \alpha^m \left[\varphi(i,j,t,t') - \sum_{v \in V}\sum_{r,r' \in R} x_{i,j,t,t',r,r'}^v\right]\right\}$$

    End

End

2.3: update $\theta_k^m$ with current iteration index and current max profit paths solution by sub-gradient method:

For each candidate recharging station $k$

    For each link $(k, j, t, t')$ which connect to recharging station $k$

$$\theta_k^m = \max\left\{0, \theta_k^{m-1} + \alpha^m \left[\sum_{v \in V}\sum_{t \in T, r \in R}\sum_{(j,t',r') \in N} x_{k,j,t,t',r,r'}^v - p_k\right]\right\}$$

    End

End

**Step 3: Solve sub-problem SPw**

With the value of multiplier $\theta_k^m$ in current iteration, solve the knapsack problem **SPw** with dynamic programing algorithm. Save the recharging station solution $w_k^m$ of **SPw** of current iteration.

**Step 4: Lower bound calculation**

4.1: With the Lagrangian multipliers $\varepsilon_{i,j,t,t'}^m$ and $\theta_i^m$ of current iteration. Update $profit^m(i,j,t,t')$ on each link according to Eq. (11)

4.2: Calculate the total profit of current iteration using $\sum_{(i,j,t,t',r,r')\in E} \varepsilon_{i,j,t,t'}^m \times \varphi(i,j,t,t')$

4.3: Call **Algorithm 2** to obtain optimal route for all vehicles.

4.4: Active virtual vehicles to serve all unserved demand.

4.5: Tally the profit of all vehicles, cost of virtual vehicles and total profit (Eq. (10)) to obtain the objective function value $z_l^m$ of **SPx** as a Lower Bound estimate of the current best solutions;

$z_l^m = \max(z_l^m, z_l^{m-1})$

**Step 5: Upper bound calculation**

5.1: Use $w_k^m$ as the recharging station solution, call **Algorithm 2** to find feasible routes for all vehicles.

5.2: Search for unsatisfied demands and use virtual vehicles to find routes to serve for all the demands which are not served by feasible routes.

5.3: Calculate the total transportation cost without multiplier values of all feasible routes.

5.3: Tally the total transportation cost including feasible routes and virtual vehicle routes to obtain the objective function value $z_u^m$ of **RRS-LRP** model 1 as the upper bound of current solutions. $z_u^m = \min(z_u^m, z_u^{m-1})$

**Step 6: Calculate the optimal gap and check Terminal condition**

Calculate $gap^m = (z_u^m - z_l^m)/z_u^m$

If $m > Q$ or $z_u^m - z_u^{m-1} < \Delta z^*$ in which Q is the maximum number of iterations.

**Algorithm terminate**

Else go to **Step 2**.

## 5. Comparison between proposed model and alternative formulations

This section aims to systematically examine different formulations for the RRS-LRP, with additional notations in Table 8 to be used in the alternative models that do not use the RST network representation**Notation**

Table 8. Additional parameter and variable definition (without using RST network)

| Symbols | Definition |
|---|---|
| $p_j^v$ | The amount of remaining resource when vehicle v arrive at node j |
| $p'^v_j$ | The amount of remaining source when vehicle v leaves node j |
| $x_{i,j}^v$ | Binary vehicle routing decision variable, =1, If vehicle v travels through link(i, j); =0, otherwise. |
| $c_i$ | Resource price at recharging station i |
| $f_v$ | Travel frequency per time period of vehicle v |
| $n_v$ | Number of stations passed by the pre-determined route of vehicle v |
| s | Subscript of the sequential index station in the pre-determined vehicle route |
| $q_{i,s}^v$ | =1 if node i is the $s^{th}$ node passed by vehicle v; =0, otherwise. |
| $w_s^v$ | Amount of resource purchased at the $s^{th}$ station on the path of vehicle v |
| $x_s^v$ | Binary variable, =1 if vehicle v stops at the $s^{th}$ station; =0, otherwise |

**(1) Basic RRS-LRP Model without time related constraints**

Worley et al. (2012) and Yang and Sun (2015) modeled the charging station location routing location routing problem with capacitated EVs, specifically for commercial EVs and battery swap stations. Their models in these two papers are developed for RRS-LRP without either time windows or space-time networks, and we further adapt them for a clear side-by-side comparison.

**Model 2.**

**Objective function**

$$z_2 = \min \sum_{v \in V} \sum_{(i,j) \in E^p} c_{i,j} x_{i,j}^v + \sum_{k \in N_r^p} b_k w_k \qquad (11)$$

**Subject to.**

$$\sum_{j \in N^p \setminus \{o\}, i \neq j} x_{i,j}^v - \sum_{j \in N^p \setminus \{d\}, i \neq j} x_{j,i}^v = 0, \forall i \in N^p \setminus \{o, d\}, \forall v \in V \qquad (12)$$

$$\sum_{j \in N^p \setminus \{o\}} x_{o,j}^v - \sum_{j \in N^p \setminus \{d\}} x_{j,d}^v = 0, \forall v \in V \qquad (13)$$

$$\sum_{v \in V} \sum_{j \in N^p \setminus \{d\}, j \neq k} x_{j,k}^v \leq p_k w_k, \forall k \in N_r^p \qquad (14)$$

$$\sum_{j \in N^p \setminus \{o\}} x_{o,j}^v \leq 1, \forall v \in V \qquad (15)$$

$$\sum_{v \in V} x_{i,j}^v \geq \varphi(i,j), \forall \varphi(i,j) \in \Psi_p \qquad (16)$$

$$\sum_{v \in V} \sum_{i,j \in S} x_{i,j}^v \leq |S| - 1, S \subset N^p, |S| > 1 \qquad (17)$$

$$p_j^v \leq {p'}_i^v - r_{i,j} x_{i,j}^v + q_v (1 - x_{i,j}^v),$$

$$\forall i \in N^p \setminus \{d\}, \forall j \in N^p \setminus \{o\}, i \neq j, v \in V \qquad (18)$$

$${p'}_o^v = q_v, \forall v \in V \qquad (19)$$

$${p'}_k^v = q_v w_k, \forall v \in V, \forall k \in N_r^p \qquad (20)$$

$$p_j^v = {p'}_j^v, \forall v \in V, \forall j \in N^p \setminus N_r^p \qquad (21)$$

$$p_j^v, {p'}_j^v \geq 0, \forall v \in V, \forall j \in N^p \qquad (22)$$

$$x_{i,j}^v \in \{0,1\}, \forall v \in V, \forall i \in N^p \setminus \{d\}, \forall j \in N^p \setminus \{o\} \qquad (23)$$

$$w_k \in \{0,1\}, \forall k \in N_r^p \qquad (24)$$

The objective function (11) minimizes the total cost including the total shipment cost and RRS construction cost. There are a number of basic constraints for the VRP such as flow balance constraints (12) and (13). constraints (14) impose the recharging station capacity constraints. Constraints (15) ensure each vehicle is assigned to at most one trip. Constraints (16) guarantee each customer must be visited by at least one vehicle. Constraints (17) are adopted to eliminate sub-tours. Among resource related constraints (18-22), Constraints (18) track the resource consuming status when vehicles arrive at nodes. Resource recharging constraints (19) and (20) represent that the vehicles are charged only at the origin and the selected recharging station. Constraints (21) are used to represent the resource flow balance for general nodes. By assisting non-negative constraints (22), the vehicles never run out of resource. Integer variables are defined by constraints (23-24).

In the above RRS-LRP model, two sets of variables are used to track the entering and exiting resource levels of a vehicle at each node. Without a time-expanded network, a large number of sub- tour elimination constraints are intended to meet the path feasibility requirement. In the real world cases, vehicles can recharge at the same station in tour if necessary; however, the formulation in model 2 prevents the vehicles from visiting the same recharging station for more than once. To allow multi-times recharging on the same RRS, the basic model needs to be further extended (Yang and Sun, 2015).

**(2) RRS-LRP model with pre-determined vehicle routes**

With all pre-determined routes, Nourbakhsh and Ouyang (2010) presented a locomotive fueling strategy optimization formulation to determine the fuel station to contract so as to minimize the sum of fuel purchasing costs, train delay cost, and contract cost. According to the model by Nourbakhsh and Ouyang, (2010), a RRS-LRP model with pre-determined routes can be stated as model 3.

**Model 3.**

**Objective function**

$$z_3 = min \sum_{v \in V} \sum_{s=1}^{n_v} f_v \left[ \sum_{i \in N^p} (c_i q_{i,s}^v w_s^v) \right] + \sum_{k \in N_r^p} b_k w_k \qquad (25)$$

**Subject to.**

$$w_s^v \leq q_v x_s^v, \forall v \in V, \forall s = 1,2,\ldots,n_v \tag{26}$$

$$g_v + \sum_{s=1}^{k-1}(w_s^v - r_{s,s+1}) + w_k^v \leq q_v, \forall v \in V, \forall k = 1,2,\ldots,n_v \tag{27}$$

$$g_v + \sum_{s=1}^{k}(w_s^v - r_{s,s+1}) \geq 0, \forall v \in V, \forall k = 1,2,\ldots,n_v - 1 \tag{28}$$

$$\sum_{s=1}^{k}(w_s^v - r_{s,s+1}) \geq 0, \forall v \in V, \forall k = n_v \tag{29}$$

$$\sum_{v \in V}\sum_{s=1}^{n_v} f_v q_{i,s}^v x_s^v \leq p_k w_k, \forall k \in N_r^p \tag{30}$$

$$w_s^v \geq 0, x_s^v \in \{0,1\}, \forall v \in V, \forall s = 1,2,\ldots,n_v \tag{31}$$

$$w_k \in \{0,1\}, \forall k \in N_r^p \tag{32}$$

The objective function (25) refers to Nourbakhsh and Ouyang, (2010) to minimize total resource cost and contract cost. The resource recharging constraints (26) ensure the vehicles stop at a selected recharging station before recharging; Constraints (27-29) are resource recharging and consumption constraints that stipulate the resource of vehicle will never be empty or exceed the capacity while traveling; constraints (30) serve as the recharging station capacity constraint the construction cost constraint is eliminated because it is included in the objective function. Constraints (31-32) are the variables' domain.

In Model 3, by assuming the predetermined routes, the flow balance constraints and sub-tour elimination constraints are needed. The resource changing status is tracked as part of the resource balance control. The construction cost can be either included in the objective function or constrained with a total budget. However, because of the predetermined route strategy, the flexibility of routes could not be fully taken into account especially under dynamic demand or time-varying congestion scenarios.

Different formulations for RRS-LRP are compared in Table 9.

Table 9. Comparison of different types of RRS-LRP formulation

| Type of constraint | Basic RRS-LRP Model 2 | Route predetermined Model 3 | Model 1 (RST network) |
|---|---|---|---|
| Flow balance | $\mathbf{Ax} = \mathbf{b}$ | predetermined Route | $\mathbf{Ax} = \mathbf{b}$ |
| Recharging station capacity | $\sum \mathbf{x} \leq \mathbf{wp}$ | $\sum \mathbf{x} \leq \mathbf{wp}$ | $\sum \mathbf{x} \leq \mathbf{wp}$ |
| Demand satisfying | $\mathbf{x} \geq \mathbf{od}$ | Not included | $\mathbf{x} \geq \mathbf{od}$ |
| Construction budget | Included in the objective function | Included in the objective function | $\sum \mathbf{bw} \leq \mathbf{B}$ |
| Subtour elimination | Yes | No | Not needed in an acyclic network |
| Resource consumption | $\mathbf{y} \leq \mathbf{y}' - \mathbf{rx} + \mathbf{q}(1 - \mathbf{x})$ <br> $\mathbf{y} > \mathbf{0}$ | $0 \leq r^+ - r^- \leq q$ | Implicitly considered in flow balance |
| Resource recharging | $\mathbf{y}' = \mathbf{q}$ | | Implicitly considered in flow balance |
| Time consumption | Not considered | Not considered | Coded in space-time network |

With the 3-dimensional network, the RRS-LRP can be formulated with only a small set of LRP constraints that allow the following rich set of features.:

(1) Travel time sensitive resource consumption. The speed-dependent nonlinear resource consumption process on the same physical link can be easily modeled by building different time-dependent resource consumption arcs.

(2) Various types of e time-dependent resource recharging links can be coded in a network to reflect different degrees of recharging efficiency within a allowable time budget.

## 6. Numerical Experiment

In this section we report 2 illustrative experiments for sensitivity analysis and 2 large scale experiments, which were conducted on an Alienware 15 personal computer running Windows 8.1 with an Intel Core i7 2.5GHz processor and 16GB of main memory. Our Lagrangian relaxation algorithm framework, dynamic programming algorithm for VRP-RS and DP algorithm for Knapsack problem are all implemented in C++ and have been compiled using the .NET framework version 4.5.51641. The source code is uploaded to GitHub website (https://github.com/GonguanLu/Lagrangian-Relaxation-Algorithm-For-RRSLRP) for references.

**Experiment 1.** To demonstrate how Lagrangian relaxation algorithm framework updates the location and vehicle routes decision during the iterative optimization, a network as Fig. 6 is used. Travel time between nodes is labeled on the link, assume that the resource cost is twice the value of travel time on each link except the recharging link between predecessor and successor recharging nodes, and the recharging link will add up to 20 units of resource to vehicle if travel through it. The depot node is the origin and destination of the vehicle. The on-link demand $\varphi(i,j,t,t')$, vehicle setting and construction cost/budget are listed in Tables 10-12, respectively.

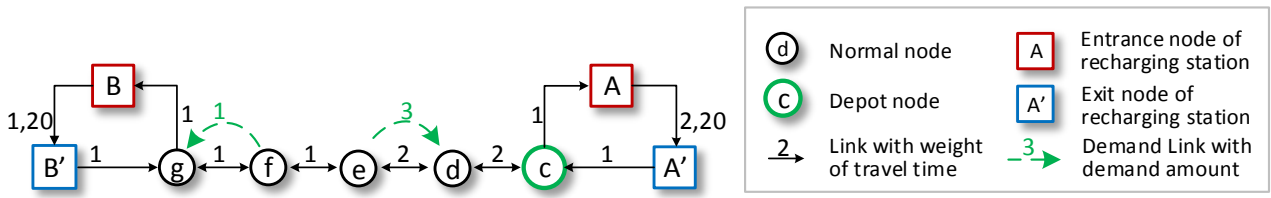

Figure 6. Test network in experiment 1

Table 10. demand links of Experiment 1

| Demand | Origin | Destination | Departure time | Arrival time |
| --- | --- | --- | --- | --- |
| 1 | f | g | 10 | 11 |
| 2 | e | d | 15 | 17 |
| 3 | e | d | 18 | 20 |
| 4 | e | d | 23 | 25 |

Table 11. Setting of vehicle of experiment 1

| Origin | Destination | Departure time window | Arrival time window | Resource capacity | Initial resource |
|---|---|---|---|---|---|
| c | c | (1,2) | (30,30) | 40 | 15 |

Table 12. Construction cost and budget of experiment 1

| Construction cost of RRS A | Construction cost of RRS B | Total construction budget |
|---|---|---|
| 10 | 10 | 15 |

Obviously, the problem is to select one out of two candidate resource recharging stations and further provide service to demand links. We show the computational data in 2 iterations, including the location decision, value of multipliers and vehicle routes of lower bound calculation, to demonstrate the iterativeprocess in Table 13, with the corresponding vehicle routes and the resource change process of lower bound in Fig. 7

Table 13. The values of multipliers in iteration 5

| Demand link/RRS | Demand amount/Capacity | Supply/Usage of last iteration | Multiplier value of current iteration |
|---|---|---|---|
| $\varphi(f,g,10,11)$ | 1 | 0 | $\varepsilon^5_{f,g,11,12} = 16.125$ |
| $\varphi(e,d,15,17)$ | 1 | 1 | $\varepsilon^5_{e,d,15,17} = 27.592$ |
| $\varphi(e,d,18,20)$ | 1 | 0 | $\varepsilon^5_{e,d,18,20} = 17.917$ |
| $\varphi(e,d,23,25)$ | 1 | 1 | $\varepsilon^5_{e,d,23,25} = 15.05$ |
| RRS A | 3 | 0 | $\theta^5_A = -7.28$ |
| RRS B | 3 | 2 | $\theta^5_B = -6.28$ |
| RRS decision | A | | |

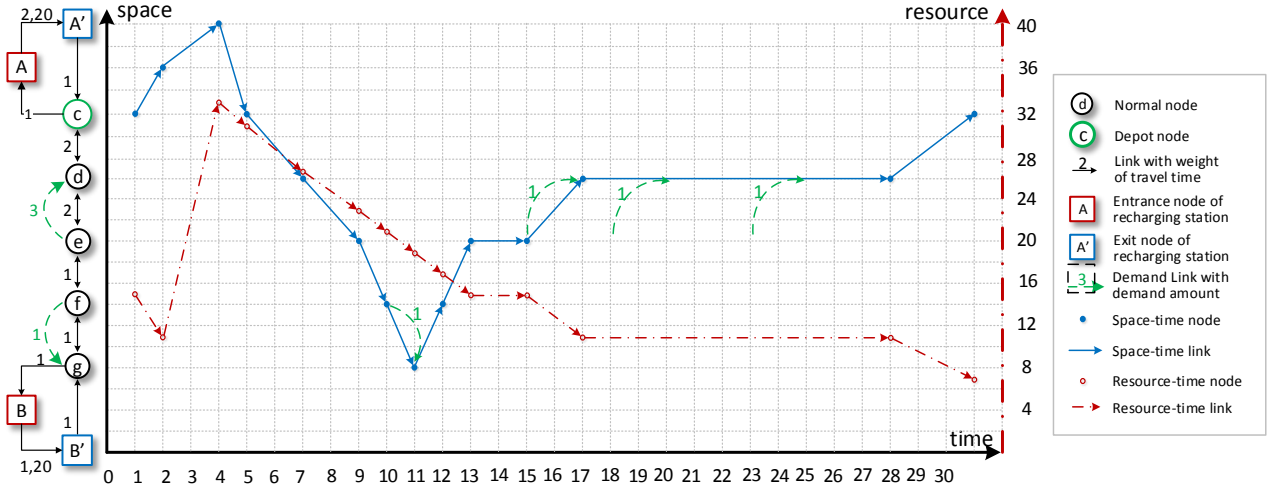

Figure 7. Routes of lower bound in iteration 5

Node A is the RRS location decision of iteration 5, and the vehicle enters RRS A for recharging and 2 demand links are served with transportation cost of 16. The dash line tracks the resource status during the route. Table 14 shows the updated multipliers are after iteration 6.

Table 14. The value of multipliers in iteration 6

| Demand link/RRS | Demand amount/Capacity | Supply/Usage of last iteration | Multiplier value of current iteration |
|---|---|---|---|
| $\varphi(f,g,10,11)$ | 1 | 1 | $\varepsilon^6_{f,g,11,12} = 16.125$ |
| $\varphi(e,d,15,17)$ | 1 | 1 | $\varepsilon^6_{e,d,15,17} = 27.592$ |
| $\varphi(e,d,18,20)$ | 1 | 0 | $\varepsilon^6_{e,d,18,20} = 21.5$ |
| $\varphi(e,d,23,25)$ | 1 | 0 | $\varepsilon^6_{e,d,23,25} = 18.633$ |
| RRS A | 3 | 1 | $\theta^6_A = -7.35$ |
| RRS B | 3 | 0 | $\theta^6_B = -7.35$ |
| RRS decision | B | | |

We can observe that, the multipliers of demands $\varphi(f,g,10,11)$ and $\varphi(e,d,15,17)$ keep no change as they have been served in the last iteration, while the other 2 unserved demand requests $\varepsilon^6_{e,d,18,20}$ and $\varepsilon^6_{e,d,23,25}$ lead to increased profit/LR multipliers.. , which lead to a selection of node B as the RRS location.

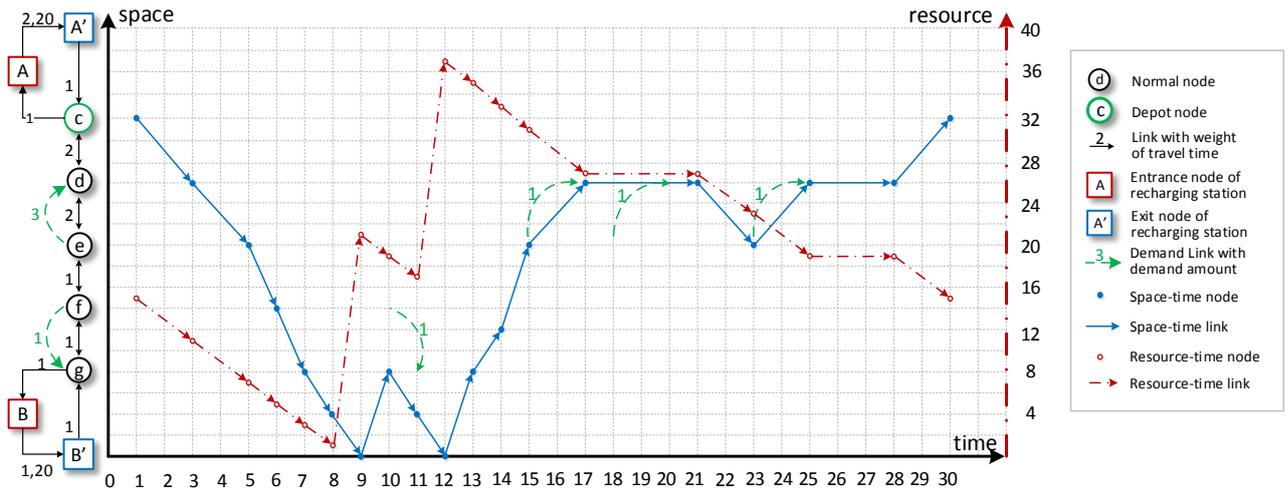

Figure 8. Routes of lower bound in iteration 6

In Fig. 8, we can see that the vehicle enters RRS B twice to obtain enough resource for traveling back to depot *c*, and 2 demand links are served by the vehicle with a total transportation cost of 22. Overall, the algorithm converges in 10 iteration and the optimal solution of experiment 1 is RRS A.

**Experiment 2.** Based on the network instance in Fig. 1, network in Fig. 9 is generated by adding resource recharging stations, recharge links and demand links.

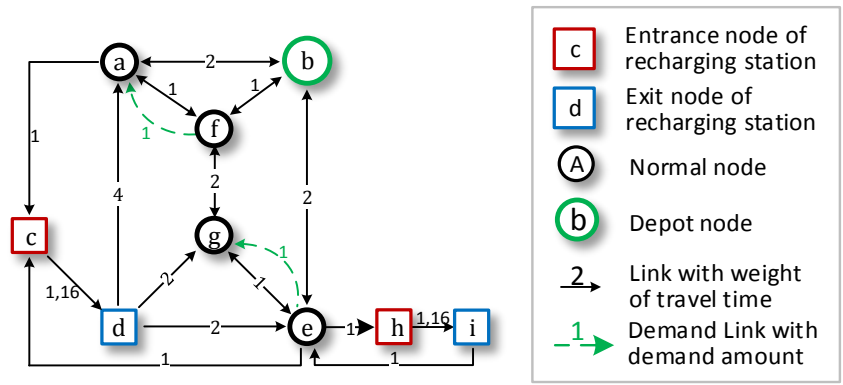

Figure 9. A network with 11 nodes, 3 RRSs and 2 demand links

Similar with experiment 1, travel time between nodes are labeled on the links, resource cost on link is set as twice as the value of travel time. Recharging links can provide at most 16 units of resource for vehicles. 2 time-dependent demands in the instance are generated in the following table. The construction cost of RRS c and h is 11 and 12, respectively, the capacity of both RRSs is set to 2. Given demand and vehicle data in Tables 15 and 16, the problem of the instance aims to determine, with the total construction budget, which RRS should be built so the system can satisfy all demand within minimum transportation cost.

Table 15. Parameter of demand in experiment 2

| Demand | Origin | Destination | Departure time | Arrival time |
|---|---|---|---|---|
| 1 | e | g | 8 | 9 |
| 2 | f | a | 3 | 4 |

Table 16. Setting of vehicles

| Origin | Destination | Departure time window | Arrival time window | Initial resource |
|---|---|---|---|---|
| B | b | (1,2) | (16,20) | 7 |

Using the vehicle set in which |V|=2, Fig. 10 shows the iteration-by-iteration converging patterns of upper bound and lower bounds. The lower bound estimates are improved smoothly after iteration 3, the relative gap can be reduced to 0 in a few iteration in this small scale experiment.

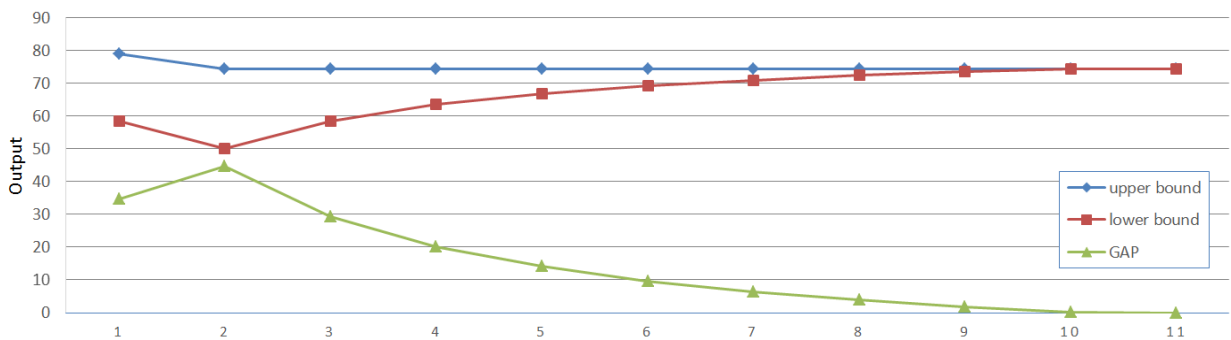

Figure 10. Evolution of upper bound and lower bound of RRS-LRP in experiment 2

**Impact of resource recharging station capacity**. With sufficient construction budget, we use 2 scenarios with different capacity of resource recharging stations as shown in Table 17 to perform a sensitivity analysis. In scenario 1, each recharging station can only serve one vehicle due to the limited capacity, thus as in the figures of vehicle routes, the 2 vehicles recharge in the different recharging stations on their way to serve for demand. For the capacity of station h is set to 0 in scenario 2, both the vehicles use the RRS c for recharging. Obviously, due to the capacity and location of RRS, vehicle 2 in scenario 1 selects RRS h for recharging so the total transportation cost (travel time) is less than in scenario 2.

Table 17. Impact of recharging station capacity

| Scenario No. | 1 | 2 |
| --- | --- | --- |
| Capacity of RRS *c* | 1 | 2 |
| Capacity of RRS *h* | 1 | 0 |
| Solution of RRS | 3,8 | 3 |
| Route of vehicle 1 | | |
| Route of vehicle 2 | | |
| value of objective function | 18 | 19 |

**Impact of construction budget**. Shown in Table 18, the RRS can be built when the budget is at least 11. In the initial setting, RRS h is more expensive than RRS c, however, even if the budget is sufficient to build RRS h, the algorithm still selects RRS d rather than h as the final decision. This behavior can be explained by that the fact that the vehicle cannot find a feasible route to serve demand (f, a) with recharging in RRS h, despite the travel time for the vehicle to serve demand (e, g) and recharge in RRS h is shorter. Once the budget increases to 25, both the RRS 3 and 8 are selected for building decision and the value of objective function drops from 19 to 18.

Table 18. Impact of construction budget

| **Construction budget** | <11 | 11 | [12, 24] | 25 |
| --- | --- | --- | --- | --- |
| Solution of RRS | null | c | c | c,h |
| Remainder budget | <11 | 0 | 1-13 | 0 |
| value of objective function | 0 | 19 | 19 | 18 |

**Impact to route solution by different resource change strategies.** Let us consider a scenario of higher fuel consumption for higher speed, where for each link whose travel time is greater than 1, vehicle can save 1 unit of time by consuming extra 2 units of resource. Table 19 lists the node, time and resource sequence of normal route without high fuel consumption strategy, and Fig. 11 shows the vehicle cannot pick up both demand once. With the optional strategy, as shown in Table 21, the vehicle can travel fast on link (*d*, *e*) and (*g*, *f*) by cost 6 units of resource, and further pick up demand requests $\varphi(e, g, 8, 9)$ and $\varphi(f, a, 3, 4)$, as shown in Fig. 12.

Table 19. Vehicle route in do nothing scenario

| Node sequence | *b* | *b* | *f* | *a* | *c* | *d* | *g* | *f* | *b* | *b* |
| --- | --- | --- | --- | --- | --- | --- | --- | --- | --- | --- |
| Time sequence | 1 | 2 | 3 | 4 | 5 | 7 | 9 | 11 | 12 | 20 |
| Resource sequence | 6 | 6 | 4 | 2 | 0 | 16 | 14 | 10 | 8 | 8 |

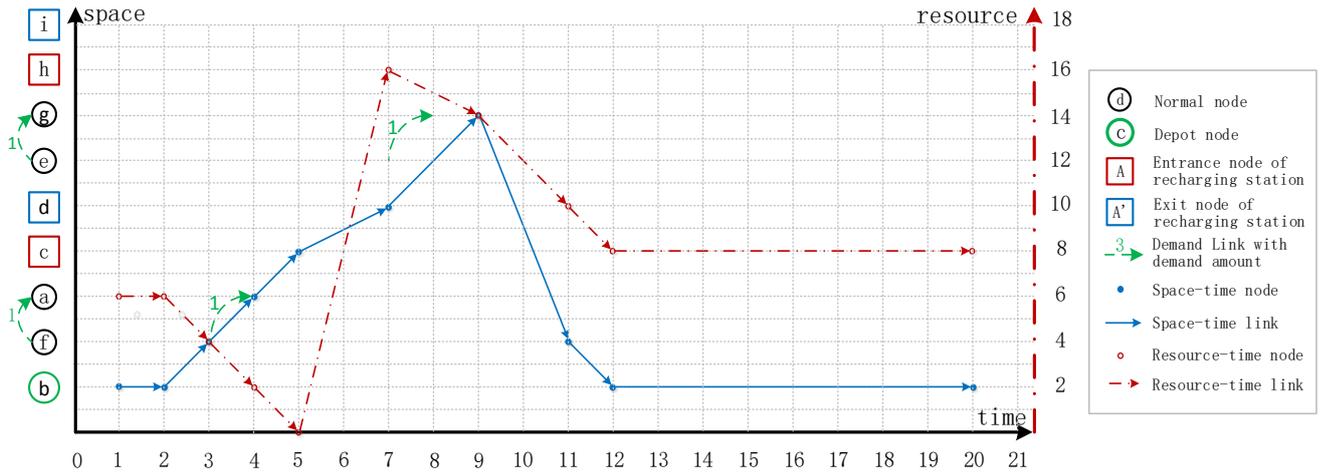

Figure 11. Vehicle routing and resource consumption process without higher fuel consumption

strategy

Table 21. Vehicle route with t with higher-fuel-consumption-for-higher-speed strategy

| Node sequence | b | b | f | a | c | d | E | g | F | b | b |
|---|---|---|---|---|---|---|---|---|---|---|---|
| Time sequence | 1 | 2 | 3 | 4 | 5 | 7 | 8 | 9 | 10 | 11 | 20 |
| Resource sequence | 6 | 6 | 4 | 2 | 0 | 16 | 10 | 8 | 2 | 0 | 0 |

*Save 1 unit of time by consuming extra 2 units of resource

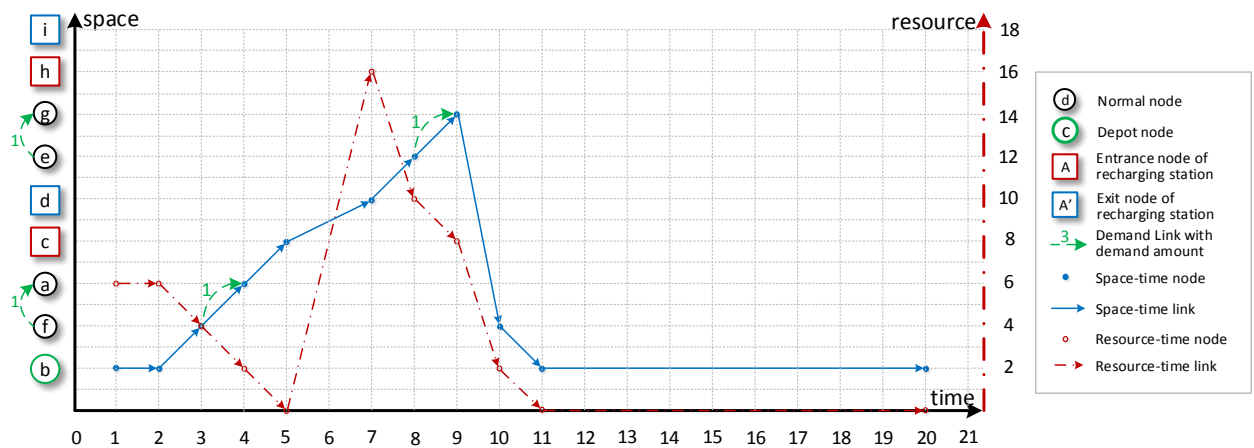

Figure 12. Vehicle routing and resource consumption process with higher-fuel-consumption-for-

higher-speed strategy

Instead of assuming a maximum level strategy for resource inventory, Table 22 shows a result for the order up-to level strategy that does not require the vehicle to refuel back to the resource capacity every time, while the amount of resource to be recharged is assumed to be proportional to the recharging time. It is interesting to observe that, in order to catch up with the time of demand link $\varphi(e, g, 8, 9)$, the vehicle opts to only recharge up to half of its resource capacity in RRS $c$, and then recharge another half capacity of resource in RRS $h$ to obtain enough resource to travel back to destination, with the

Table 22. Vehicle route with order-up-to-level strategy

| Node sequence | b | b | F | a | c | d | e | g | C | d | G | f | b | b |
|---|---|---|---|---|---|---|---|---|---|---|---|---|---|---|
| Time sequence | 1 | 2 | 3 | 4 | 5 | 6 | 8 | 9 | 11 | 12 | 14 | 16 | 17 | 20 |
| Resource sequence | 6 | 6 | 4 | 2 | 0 | 8 | 4 | 2 | 0 | 8 | 6 | 2 | 0 | 0 |

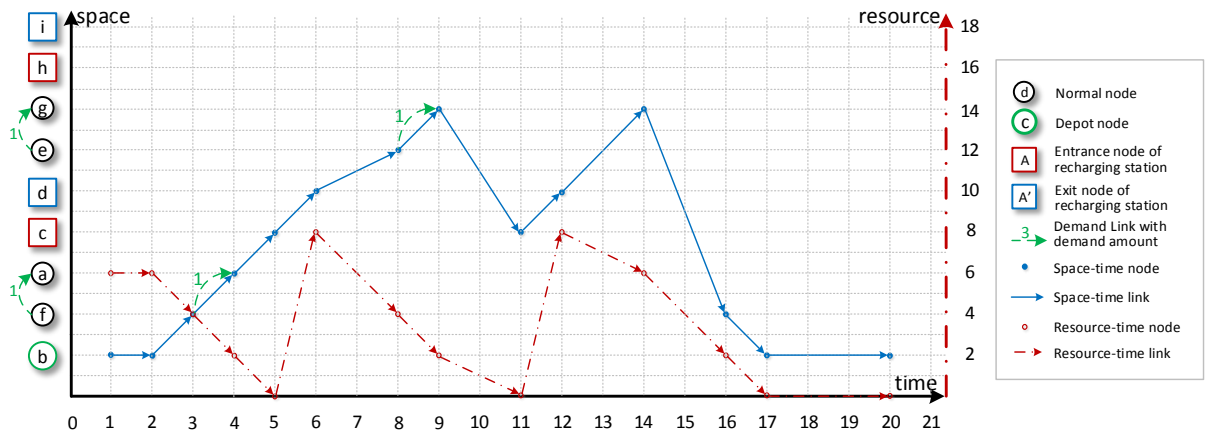

Figure 13. Vehicle routing and resource consumption process with order-up-to-level strategy

**Experiment 3.** We now consider a simplified Sioux Falls network shown in Fig. 14, which consists of 29 nodes and 81 links. Specifically, there are 12 demand links, 15 vehicles are assigned on nodes with initial resource of 15 units. We have 5 candidate RRSs with average construction cost of 20 and capacity of 5, with a total construction budget B = 60.

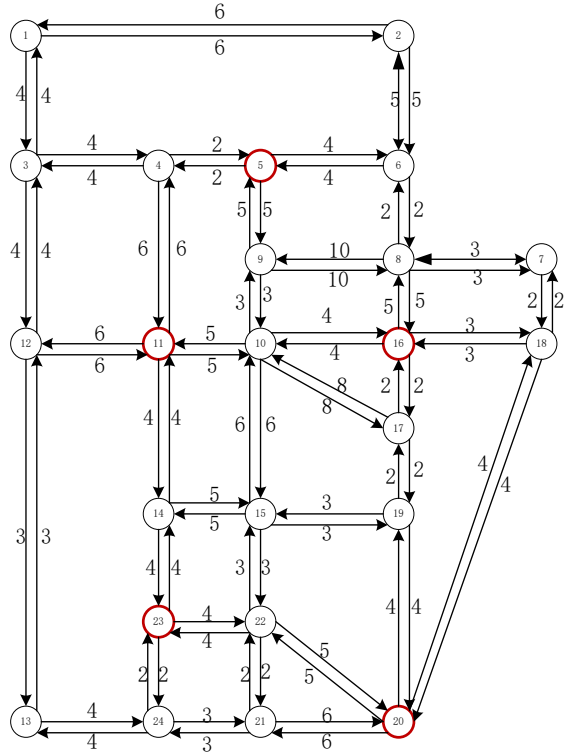

Figure 14. Sioux Falls network with 5 candidate RRSs

Our implemented algorithm takes 88.26 seconds, approximately 1.77 second per iteration to find a reasonably good RRS location solution within a 13.1% gap. The final depot solution is a set of nodes 11, 16 and 23. Intuitively, the nodes in the center of the network offers better opportunity for vehicles to recharge within less traveling distance. The solution gap pattern is reported in Fig. 15, and the dataset and problem solution of this RRS-LRP experiment is available at our GitHub project [web site](web site). The significant gap is mainly caused by many possibly overserved demand requests, that is, the relaxed demand satisfaction constraints still allow multiple vehicles to pick up a single high-profit demand for more than once, due to its profit-maximization nature in the lower bound routing solutions. A branch and bound solution method that assigns different vehicles to service requests precisely once could be beneficial to further reduce the gap.

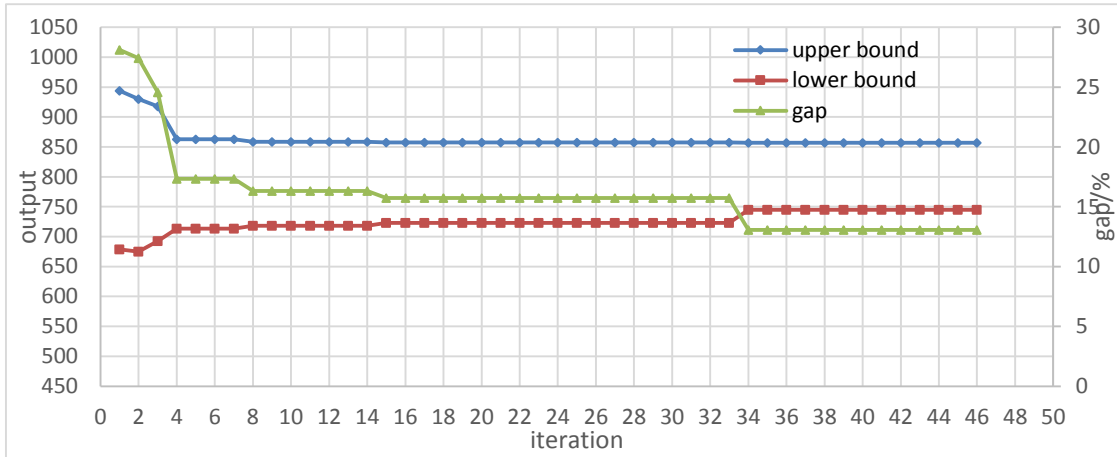

Figure 15. Evolution of upper bound and lower bound of RRS-LRP in experiment 3

**Experiment 4.** To test the efficiency of the Lagrangian Relaxation Algorithm Framework, A real world network case of Chicago is tested (http://www.bgu.ac.il/~bargera/tntp/) with 933 nodes and 2967 links shown in Fig. 16. Some inputs are:

- 40 demand links

- 30 vehicles are assigned on nodes with initial resource of 80.

- 40 candidate RRSs (red dot in Figure 16) with construction cost around 10 and service capacity of 3

- Construction budget B = 200.

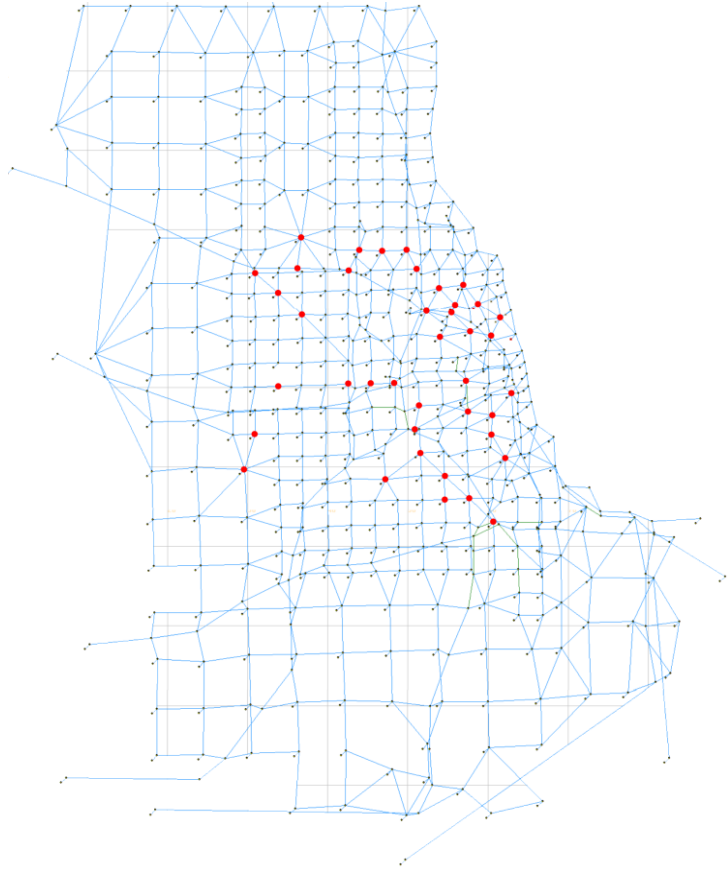

Figure 16. The Chicago network with 40 hypothetic candidate RRSs

This network is an aggregated representation of the Chicago region, and our hypothetical case selects about 20 RRSs out of 40 candidate recharging stations. In a single CPU thread, the algorithm takes 87.2 minutes of CPU time, approximately 1.34 minutes per iteration, to find a reasonable good RRS location and vehicle routs solution with a 7.6% relative gap, as reported in Fig. 17. The average calculating time for a single vehicle routing process is 3.18 seconds. The efficiency of our proposed algorithm is much better than what we implement in GAMS with the commercial solver, and the CUP run time can be further reduced by applying parallel computing with multi-core CPU for the decomposed routing algorithms for different vehicles.

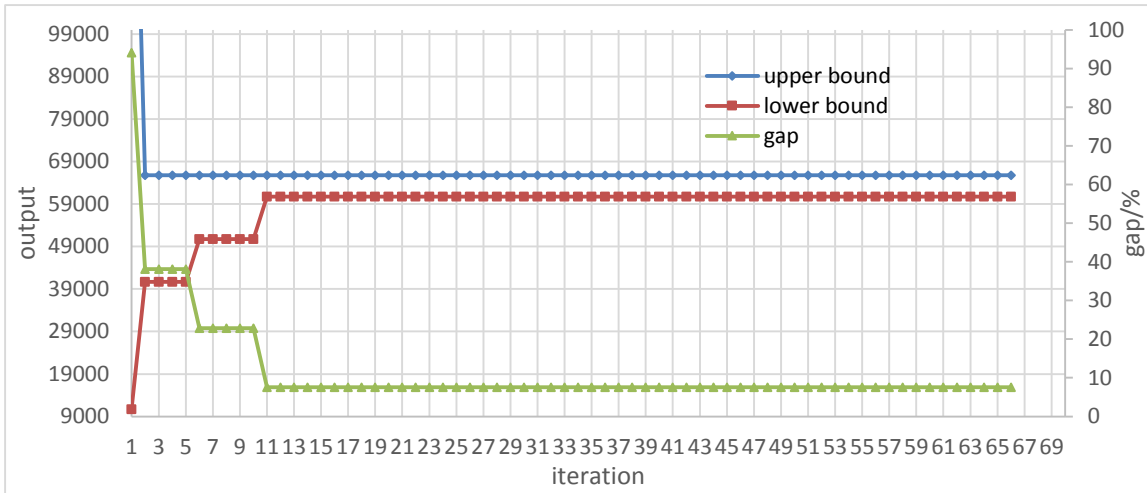

Figure 17. Evolution of upper bound and lower bound of RRS-LRP in experiment 4

## 7. Conclusion

Motivated by emerging research needs of better designing resource recharging service network we studied the resource recharging station location routing problem. By introducing a directed acyclic resource-space-time network, we propose a multi-commodity flow model for RRS-LRP. Many complex constraints, such as resource consumption/charging, time window restrictions for demand requests are directly coded in our three-dimensional well-structured network representation. We also compare our formulation with existing literature, and highlight the potential for considering a rich set of many practically important features in vehicle routing and resource recharging/consumption activities.

We developed a Lagrangian relaxation algorithm framework to decompose the original problem into knapsack subproblem and vehicle routing problems, which are solvable through two dynamic programming algorithms. We systematically test the developed algorithms for different networks. In our future research, we need to consider more generic demand representations to model the pickup-delivery requirements (with possible different commodities), instead of the single-purpose demand link. Secondly, a more comprehensive inventory location-routing modeling framework in the RST network could represent real-world vehicle carrying capacity and warehouse inventory capacity as particular types of resources. Thirdly, we need to either use heuristic algorithms to find better upper bound solutions or embed a branch and bound algorithm to better enforce the demand-to-vehicle assignment constraints, so as to reduce the solution gaps for large-scale test cases.

**Acknowledgements**

The research reported here was partially sponsored by the National United Engineering Laboratory of Integrated and Intelligent Transportation, Southwest Jiaotong University, China. The second author is partially supported by National Science Foundation – United States Grant No. CMMI **1538105** "Collaborative Research: Improving Spatial Observability of Dynamic Traffic Systems through Active Mobile Sensor Networks and Crowdsourced Data", and project RCS2015K006 sponsored by State Key Laboratory of Rail Traffic Control and Safety, Beijing Jiaotong University, China.

**Appendix A**

*Algorithm for setting up the cost matrix elements in resource-space-time network*

**Input:**

Travel time $TT_{i,j}$, travel cost $c_{i,j}$ and resource cost $r_{i,j}$ on each links in physical network

**Output:**

Travel cost $c_{i,j,t,t',r,r'}$ on each links in resource-space-time network

**Resource-space-time network generating algorithm**

**Step 1: Initialization.**

For all RST links, $c_{i,j,t,t',r,r'} = M$.

**Step 2: Weight of normal RST link cost.**

For all physical link $(i,j) \in E^p$

    For (t=0 to |T|)

        For (r=$r_0$ to |R|)

$$c_{i,j,t,t+TT_{i,j},r,r+r_{i,j}} = c_{i,j}$$

        End For

    End For

End For

**Step 3: Weight of RST waiting links.**

For all physical node $i \in N^p$

    For (t=0 to |T|)

        For (r=$r_0$ to |R|)

$$c_{i,j,t,t+1,r,r} = 0$$

        End For

    End For

End For

**Step 4: Weight of resource exhausting links.**

For all physical node $i \in N^p$

    For (t=0 to |T|)

        For (r=$r_0$ to |R|)

            $c_{i,j,t,t,r,r_0} = 0$

        End For

    End For

End For

**Algorithm End**